\documentclass{elsart}

\usepackage{amssymb}
\usepackage{amsmath}
\numberwithin{equation}{section}

\begin{document}

\begin{frontmatter}



\title{Inverse Problem For Dirichlet Boundary Value Problems }


\author{Alp Arslan K\i ra\c{c},\quad Fatma Y\i lmaz}
\ead{aakirac@pau.edu.tr}
\address{Department of Mathematics, Faculty of Arts and Sciences, Pamukkale
University, 20070, Denizli, Turkey}

\begin{abstract}
In this article we consider Sturm-Liouville operator with $q\in W_{1}^{2}[0,1]$ and Dirichlet boundary conditions. We prove that if the set $\{(n\pi)^{2}:n\in \mathbb{N}\}$ is a subset of the spectrum of the Sturm-Liouville operator with Dirichlet boundary conditions, then $q=0$ a.e.
\end{abstract}

\end{frontmatter}

\section{Introduction}
Denote by $L(q)$ the operator generated in $L_{2}(0,1)$ by the expression
\begin{equation}  \label{1}
-y^{\prime\prime}+q(x)y,
\end{equation}
and by the boundary conditions 
\begin{equation} \label{per}
y(1)=y(0)=0
\end{equation}
where $q(x)\in\,W_{1}^{2}[0,1]$ is a real-valued function and  $q^{(k)}(0)=q^{(k)}(1)$, $k=0,1$.\\ 
In \cite{Yilmaz-Veliev}, they obtained asymptotic formulas of arbitrary order for eigenfunctions and eigenvalues of the Dirichlet boundary value problem (\ref{1}) with (\ref{per}) when $q(x)$ is a complex-valued summable function.\\
In this point, we examine Ambarzumyan theorem in \cite{Ambarzumian}. In 1929, Ambarzumyan  \cite{Ambarzumian} proved the following theorem as the first theorem in inverse spectral theory:
\begin{thm}
 If $\{n^{2}:n=0,1,...\}$ is the spectrum of the Neumann boundary condition with Sturm-Liouville operator, then $q=0$ a.e.
\end{thm}
In the P{\"{o}}schel-Trubowitz \cite{Poschel-Trubowitz} works, the $\sigma=\{n^{2}:n\in\mathbb{N}\}$ spectrum for the Dirichlet problem corresponds to zero and they showed that there are infinitely many $L^{2}$ potentials near zero. For the Dirichlet problem, when the spectrum is zero, the potential does not have to be zero. Therefore, the Ambarzumyan's theorem is not valid (see \cite{Chern-Law-Wang}).\\
In \cite{Chern-Law-Wang}, putting an additional condition on the potential they expanded the classical Ambarzumyan's theorem for the Sturm-Liouville equation to the general separated boundary conditions. They obtained the following this theorem (see \cite[Theorem 1.1]{Chern-Law-Wang}).\\
\begin{thm}
The Sturm-Liouville problem
\[
-y^{\prime\prime}+qy=\lambda \,y
\]
such that 
\[
y(0)\cos \alpha+y^{\prime}(0)\sin \alpha=0
\]
\[
y(\pi)\cos \beta+y^{\prime}(\pi)\sin \beta=0
\]
where $q\in\,L^{1}(0,\pi), \alpha,\beta\in[0,\pi)$.
For the Sturm-Liouville problem, assume that $\alpha=\beta\neq\frac{\pi}{2}$. Then $\sigma=\{n^{2}:n\in\mathbb{N}\}$ and the potential function $q$ satisfies
\[
\int_{0}^{\pi}q(x)\,\cos 2(x-\alpha)\,dx=0
\]
if and only if $q=0$ a.e.
\end{thm}
In \cite{Yurko}, Yurko also proved the following generalization theorem of the Ambarzumyan theorem on wide classes of the arbitrary self-adjoint boundary conditions and self-adjoint differential operators (see \cite[Theorem 3]{Yurko}).
\begin{thm}
Let 
\[
\lambda_{0}=\tilde{\lambda}_{0}+\frac{(\hat{q}\tilde{y}_{0},\tilde{y}_{0})}{(\tilde{y}_{0},\tilde{y}_{0})},
\]
where $\tilde{y}_{0}(x)$ is an eigenfunctions $\tilde{L}$ related to $\tilde{\lambda}_{0}$. Then 
\[
q(x)=\tilde{q}(x)+\lambda_{0}-\tilde{\lambda}_{0} \quad \textrm{a.e. on} \quad (0,1).
\]
\end{thm}
K\i ra\c{c} \cite{quasi-periodic} proved for the Sturm-Liouville operators with quasi-periodic boundary conditions without putting any condition on $q$ the potential and $q$ can be determined from a single spectrum. Thus he get the classical Ambarzumyan's theorem. This result is the following (see \cite[Theorem 1]{quasi-periodic}).
\begin{thm}
 If first eigenvalue of the operator $L_{t}(q)$ for any fixed number $t$ in $[0,2\pi)$ is not less than the value of $\min\{t^{2},(2\pi-t)^{2}\}$ and the spectrum $S(L_{t}(q))$ contains the set $\{(2n\pi-t)^{2}:n\in\mathbb{N}\}$, then $q=0$ a.e.  
\end{thm}

K\i ra\c{c} \cite[Theorem 1.1]{Kirac-inverse} proved the following inverse spectral result for periodic and anti-periodic boundary conditions:
\begin{thm}
Denote the n th instability interval by $\ell_{n}$, and suppose that $\ell_{n}=o(n^{-2})$ as  $n\rightarrow\infty$. Then the following two assertions hold:
\\(i) If $\{(n\pi)^{2}: n$ even and $n>n_{0}\}$ is a subset of the periodic spectrum of the Hill operator then $q=0$ a.e. on $(0,1)$,\\
(ii) If $\{(n\pi)^{2}: n$ odd and $n>n_{0}\}$ is a subset of the anti-periodic spectrum of the Hill operator then $q=0$ a.e. on $(0,1)$.
\\Given $\epsilon>0$, there exists $n_{0}=n_{0}(\epsilon)\in\mathbb{N}$, a sufficiently large positive integer such that 
\[
\ell_{n}<\epsilon n^{-2} \quad\textrm{for all} \quad n>n_{0}(\epsilon).
\]
\end{thm}
In this article, we prove the following inverse spectral result:

\begin{thm}\label{maintheorem}
Let $q(x)\in W_{1}^{2}[0,1]$ and $q^{(k)}(0)=q^{(k)}(1)$, $k=0,1$ hold. Then the following assertion satisfies:
 
If $\{(n\pi)^{2}:n\in \mathbb{N}\}$ is a subset of the spectrum of the operator (\ref{1})-(\ref{per}) then $q=0$ a.e. on (0,1).
\end{thm}

\section{Preliminaries}
We shall consider the problem (\ref{1})-(\ref{per}). Considering Theorem 1 in \cite{Yilmaz-Veliev}, the eigenvalues of the Sturm-Liouville operator for $m \geqq N$ without loss of generality  assumption $c_{0}=0$ and 
\begin{equation}\label{cm}
 c_{m}=\int_{0}^{1}q(x) \cos m \pi x\,dx
\end{equation}
such that we obtain asymptotic formula for $\{\lambda_{m}\}$  eigenvalues
\begin{equation}\label{mc0}
\lambda_{m}=(m\pi)^{2}+c_{0}-c_{2m}+O\left(  \frac{\ln |m|}{m}\right) 
\end{equation}
form taking into account that satisfying
\begin{equation}\label{m1}
\lambda_{m}=(m\pi)^{2}+o(1).
\end{equation}
Here we denote by $N$ a large positive integer. Using this formula, for all $k\neq m$, $k=0,1,...$ the following inequality holds  
\begin{equation}\label{m2}
|\lambda_{m}-(\pi k)^{2}|>|(m-k)\pi||(m+k)\pi|-c_{1}m^{1/2}>c_{2}m
\end{equation}
for $m \geqq N$, where we denote by $c_{n}$, $n=1,2...$ , positive constants whose exact value is not essential.
To obtain the asymptotic formula for eigenvalues $\lambda_{m}$ corresponding to the normalized eigenfunctions $\Psi_{m}(x)$ of $L(q)$ we use (\ref{m2}) and the following relation 
\begin{equation}\label{m3}
(\lambda_{N}-(\pi m)^{2})(\Psi_{N}(x),\sin m\pi x)=(q(x)\Psi_{N}(x), \sin m\pi x).
\end{equation}
Here
\begin{equation}\label{m4}
\Psi_{m}(x)=\sqrt{2}\sin m\pi x+O(\frac{1}{m}).
\end{equation}

\begin{equation}\label{f5}
|(q(x)\Psi_{N}(x), \sin m\pi x)|<4M
\end{equation}
satisfy for $\forall m, \forall N\gg 1$, where $M=\int_{0}^{1}|q(x)|\,dx$.

Now, the expansion 
\[
\Psi_{N}(x)=\sum_{m_{1}>-m}^{\infty}2(\Psi_{N}(x), \sin (m+m_{1})\pi x)sin (m+m_{1})\pi x
\]
of $\Psi_{N}(x)$ by the orthonormal basis $\{\sqrt{2}\sin (m+m_{1})\pi x:m_{1}>-m\}$ has the formula 
\[
\Psi_{N}(x)=\sum_{m_{1}>-m}^{n}2(\Psi_{N}(x), \sin (m+m_{1})\pi x)sin (m+m_{1})\pi x+g(x)
\]
where $\sup_{x\in[0,1]}|g(x)|<\frac{c_{3}}{n}$.\\
 This formula is substituted in the $(q(x)\Psi_{N}(x), \sin m\pi x)$ expression and as $n\rightarrow \infty$ we obtain 
\[
(q(x)\Psi_{N}(x), \sin m\pi x)=
\]
\begin{equation}\label{m5}
\sum_{m_{1}>-m}^{\infty}2(q(x), (\sin (m+m_{1})\pi x)(\sin m \pi x))(\Psi_{N}(x),\sin(m+m_{1})\pi x).
\end{equation}
By doing the necessary operations in (\ref{m5}), we obtain the following equation (see \cite[p.155]{Yilmaz-Veliev})
\[\!\!\!\!\!\!\!\!\!\!\!\!\!\!\!\!\!\!\!\!\!\!\!\!\!\!\!\!\!\!\!\!\!\!\!\!\!\!\!\!\!\!\!\!\!\!\!(q(x)\Psi_{N}(x), \sin m\pi x)=\sum_{m_{1}=1}^{\infty}c_{m_{1}}(\Psi_{N}(x), \sin (m+m_{1})\pi x)
\]
\[
+\sum_{m_{1}=1}^{\infty}c_{m_{1}}(\Psi_{N}(x), \sin (m-m_{1})\pi x).
\]
Substituting this equality into (\ref{m3}) we get
\begin{equation}\label{m15}
(\lambda_{m}-(\pi m)^{2})(\Psi_{N}(x), \sin m\pi x)=\sum_{m_{1}=-\infty}^{\infty}c_{m_{1}}(\Psi_{N}(x), sin(m+m_{1})\pi x)
\end{equation}
where $c_{m}=\int_{0}^{1}q(x)\cos m\pi x\,dx$. Also, $c_{m}=c_{-m}$, $c_{m}\rightarrow0$ as $|m|\rightarrow\infty$.\\
Now we isolate the terms that contain the expression $(\Psi_{N}(x),\sin(m+m_{1})\pi x)$ to the right of (\ref{m15}). For this we have $m$ and $m+m_{1}$ instead of $N$ and $m$ in (\ref{m15}). Hence we obtain the following equation 
\[
\!\!\!\!\!\!\!\!\!(\lambda_{m}-(\pi m)^{2})(\Psi_{m}(x),\sin m \pi x)=-c_{2m}(\Psi_{m}(x),\sin m \pi x)
\]
\begin{equation}\label{m16}
+\sum_{m_{1},m_{2}=-\infty \atop m_{1}\neq-2m}^{\infty}\frac{c_{m_{1}}c_{m_{2}}(\Psi_{m}(x),\sin(m+m_{1}+m_{2})\pi x)}{\lambda_{m}-(\pi(m+m_{1}))^{2}}.
\end{equation}
Again, we isolate the terms for which $m_{1}+m_{2}=0,-2m$. First, $-m_{1}$ is written instead of $m_{2}$ and use the equality $c_{-m_{1}}=c_{m_{1}}$ and the last sum of this equation is isolated for which $m_{1},m_{2},m_{3}=0,-2m$ (see (12), (13) of \cite{Yilmaz-Veliev}) we get the following lemma for our notation. 
\begin{lem}\label{L1}(see \cite{Yilmaz-Veliev})
The eigenvalue $\lambda_{m}$ of the operator $L(q)$
satisfies the asymptotic formula 
\begin{equation}\label{m6}
\lambda_{m}=(m\pi)^{2}+c_{0}-c_{2m}+a_{1}(\lambda_{m})-b_{1}(\lambda_{m})+a_{2}(\lambda_{m})-b_{2}(\lambda_{m})+R_{3} 
\end{equation}
where $q(x)$ is a real-valued summable function and $c_{m}=\int_{0}^{1}q(x)\,\cos m\pi x\,dx$,
\begin{equation}\label{e6}
a_{1}(\lambda_{m})=\sum_{m_{1}=-\infty \atop m_{1}\neq0,-2m}^{\infty}\frac{c_{m_{1}}c_{m_{1}}}{[\lambda_{m}-(\pi(m+m_{1}))^{2}]}
\end{equation}
\begin{equation}\label{e8}
b_{1}(\lambda_{m})=\sum_{m_{1}=-\infty \atop m_{1}\neq0,-2m}^{\infty}\frac{c_{m_{1}}c_{2m+m_{1}}}{[\lambda_{m}-(\pi(m+m_{1}))^{2}]}
\end{equation}
\begin{equation}\label{e7}
a_{2}(\lambda_{m})=\sum_{ m_{1},m_{2}=-\infty  \atop m_{1},m_{1}+m_{2}\neq0,-2m}^{\infty}\frac{c_{m_{1}}c_{m_{2}}c_{m_{1}+m_{2}}}{\prod_{t=1,2}[\lambda_{m}-(\pi(m+m_{t}))^{2}]}
\end{equation}
\begin{equation}\label{e9}
b_{2}(\lambda_{m})=\sum_{m_{1},m_{2}=-\infty \atop m_{1},m_{1}+m_{2}\neq0,-2m}^{\infty}\frac{c_{m_{1}}c_{m_{2}}c_{2m+m_{1}+m_{2}}}{\prod_{t=1,2}[\lambda_{m}-(\pi(m+m_{t}))^{2}]}
\end{equation}
\begin{equation}\label{m9}
R_{3}=\sum_{m_{1},m_{2},m_{3}=-\infty}^{\infty}\frac{c_{m_{1}}c_{m_{2}}c_{m_{3}}(q(x)\Psi_{m}(x), \sin(m+m_{1}+m_{2}+m_{3})\pi x)}{\prod_{t=1,2,3}[\lambda_{m}-(\pi(m+m_{t}))^{2}]}.
\end{equation}
\end{lem}
Here, using the following equality 
\[
\frac{1}{-m_{1}(2m+m_{1})}=\frac{1}{2m}\left( \frac{1}{2m+m_{1}}-\frac{1}{m_{1}}\right) 
\]
we get the relation 
\begin{equation}\label{m10}
\sum_{m_{1}\neq0,-2m}\frac{1}{|-m_{1}(2m+m_{1})|}=O\left( \frac{ln|m|}{m}\right).
\end{equation}
In \cite{Yilmaz-Veliev}, obtained the following estimate (see (19) of \cite{Yilmaz-Veliev}, \cite{Veliev-Duman})
\begin{equation}\label{m13}
R_{3}=O\left( \left( \frac{ln|m|}{m}\right) ^{3}\right) .
\end{equation}

\section{Main Results}\label{results}

In our next notations we will need the following relation.\\
Let us see that   
\begin{equation}\label{f1}
c_{2m_{1}}=\int_{0}^{1}\tilde{q}(x)\,e^{-i(2m_{1})\pi x}\,dx
\end{equation}
and
\begin{equation}\label{f2}
c_{2m_{1}+1}=\int_{0}^{1}\hat{q}(x)\,e^{-i(2m_{1}+1)\pi x}\,dx,
\end{equation}
equalities are satisfying by using the even and odd functions 
\begin{equation}\label{f3}
 \tilde{q}(x)=\frac{q(x)+q(1-x)}{2},
\end{equation}
\begin{equation}\label{f4}
\hat{q}(x)=\frac{q(x)-q(1-x)}{2}
\end{equation}
in (\ref{cm}), respectively of the $q$ function in $[0,1]$.
Really, 
\[
\!\!\!\!\!\!\!\!\!\!\!\!\!\!\!\!\!\!\!\!\!\!\!\!\!\!\!\!\!\!\!\!\!\!\!\!\!\!\!\!\!\!\!\!\!\!\!\!\!\!\!\!\!\!
\int_{0}^{1}\tilde{q}(x)\,e^{-i(2m_{1})\pi x}\,dx=\int_{0}^{1}\frac{q(x)+q(1-x)}{2}\cos 2m_{1}\pi x\,dx
\]
\begin{equation}\label{d0}
-i\int_{0}^{1}\frac{q(x)+q(1-x)}{2}\sin 2m_{1}\pi x\,dx.
\end{equation}
By using $1-x=u$ change of variable in both integrals on the right-hand side of the equation (\ref{d0}), we get (\ref{f1}).
Similarly, we obtain (\ref{f2}).

\begin{lem}\label{L2}
Let $q(x)\in W_{1}^{1}[0,1]$, $q(0)=q(1)$ and $c_{0}=0$.\\
 The following asymptotic estimates are valid for the relation in (\ref{e8}) and (\ref{e9}):
\[
b_{1}(\lambda_{m})=o(m^{-2}),
\]
\[
b_{2}(\lambda_{m})=o(m^{-2}).
\]
\end{lem}
\begin{pf}
Now, for the proof of lemma by using (\ref{m3}), (\ref{f5}) and (\ref{m2}), we obtain the following asymptotic estimate (see \cite{Veliev-Duman})
\[
\!\!\!\!\!\!\!\!\!\!\!\!\!\!\!\!\!\!\!\!\!\!\! \sum_{m \in \mathbb{Z}}|(\psi_{m}(x), \sin m \pi x)|^{2}< \sum_{m \in \mathbb{Z}}\frac{(4M)^{2}}{|-m_{1}(2m+m_{1}) \pi^{2}|^{2}}
\]
\begin{equation}\label{f6}
=O\left( \frac{1}{m^{2}}\right).
\end{equation}

Arguing as in \cite[Theorem 1.2]{Kirac-inverse} arbitrary a constant $C$ with by using (\ref{m1}), (\ref{m2}) and (\ref{f6}), we obtain that
\[
\!\!\!\!\!\!\!\!\!\!\!\!\!\!\!\!\!\!\!\!\!\!\!\!\sum_{m_{1}\neq 0,-2m}\left|\frac{1}{\lambda_{m}-(m+ m_{1})^{2}\pi^{2}}-\frac{1}{m^{2}\pi^{2}-(m+ m_{1})^{2}\pi^{2}}\right|
\]
\begin{equation}\label{dif}
\leq C \quad o(1) \sum_{m_{1}\neq 0,-2m}|m_{1}|^{-2}|2m+ m_{1}|^{-2}=o(m^{-2}).
\end{equation}

First, we prove $b_{1}(\lambda_{m})$. For this, we write taking into account that (\ref{dif}) and the sum on the right-hand-side of equation (\ref{e8}) by using equalities (\ref{f1}) and (\ref{f2}), as following (see \cite[Lemma 3]{Kirac-Titchmarsh},  \cite{Shkalikov-Veliev,Kirac-inverse}).
\[
b_{1}(\lambda_{m})=\sum_{m_{1}=-\infty \atop m_{1}\neq0,-2m}^{\infty}\frac{c_{m_{1}}c_{2m+m_{1}}}{[\lambda_{m}-(\pi(m+m_{1}))^{2}]},
\]
by using (\ref{m1}) and with together equality
\begin{equation}\label{f7}
\lambda_{m}-(\pi (m+m_{1}))^{2}=-m_{1}(2m+m_{1}) \pi^{2}+o(1),
\end{equation}
\[
b_{1}(\lambda_{m})=\frac{1}{\pi^{2}}\sum_{m_{1}=-\infty \atop m_{1}\neq0,-2m}^{\infty}\frac{c_{m_{1}}c_{2m+m_{1}}}{-m_{1}(2m+m_{1})}+o(m^{-2})
\]
\[
=\frac{1}{\pi^{2}}\sum_{2m_{1}\neq 0,-2m}\frac{c_{2m_{1}}c_{2m+2m_{1}}}{-2m_{1}(2m+2m_{1})}+\frac{1}{\pi^{2}}\sum_{2m_{1}+1\neq -2m}\frac{c_{2m_{1}+1}c_{2m+(2m_{1}+1)}}{-(2m_{1}+1)(2m+2m_{1}+1)}
\]
\[
\!\!\!\!\!\!\!\!\!\!\!\!\!\!\!\!\!\!\!\!\!\!\!\!\!\!\!\!\!\!\!\!\!\!\!\!\!\!\!\!\!\!\!\!\!\!\!\!\!\!\!\!\!\!\!\!\!\!\!\!\!\!\!\!\!\!\!\!\!\!\!\!\!\!\!\!\!\!\!\!\!\!\!\!\!\!\!\!\!\!\!\!\!\!\!\!\!\!\!\!\!\!\!\!\!\!\!\!\!\!\!\!\!\!\!\!\!\!\!\!\!\!\!\!\!\!\!\!\!\!\!\!\!\!\!\!\!\!\!\!\!\!\!\!\!\!\!\!\!\!\!\!\!\!\!\!\!\!\!\!\!\!+o(m^{-2})
\]
\[
\!\!\!\!\!\!\!\!\!\!\!\!\!\!\!\!\!=-\int_{0}^{1}(\tilde{Q}(x)-\tilde{Q}_{0})^{2}\,e^{i2m\pi x}\,dx-\int_{0}^{1}(\hat{Q}(x)-\hat{Q}_{0})^{2}\,e^{i2m\pi x}\,dx+o(m^{-2})
\]
\[
\!\!\!\!\!\!\!\!\!\!\!\!\!\!\!\!\!\!\!\!\!\!\!\!\!\!\!\!\!\!\!\!\!\!\!\!\!\!\!\!\!\!\!\!\!\!\!\!\!\!\!\!\!\!\!\!\!\!\!\!\!\!\!\!\!\!\!\!\!\!\!\!\!\!\!\!\!\!\!\!\!\!\!\!\!\!\!\!\!\!\!=\frac{1}{i\pi 2m}\int_{0}^{1}2(\tilde{Q}(x)-\tilde{Q}_{0})\,\tilde{q}(x)\,e^{i2m\pi x}\,dx
\]
\begin{equation}\label{d2}
\!\!\!\!\!\!\!\!\!\!\!\!\!\!\!\!\!\!\!\!\!\!\!\!\!\!\!\!\!\!\!\!\!\!\!\!\!\!\!\!\!\!\!\!\!\!\!\!\!\!\!\!\!\!\!+\frac{1}{i\pi2m}\int_{0}^{1}2(\hat{Q}(x)-\hat{Q}_{0})\,\hat{q}(x)\,e^{i2m\pi x}\,dx+o(m^{-2})
\end{equation}
where 
\[
\tilde{Q}(x)-\tilde{Q}_{0}=\sum_{2m_{1}\neq0}\tilde{Q}_{m_{1}}\,e^{i(2m_{1})\pi x},\qquad\hat{Q}(x)-\hat{Q}_{0}=\sum \hat{Q}_{m_{1}}\,e^{i(2m_{1}+1)\pi x}
\]
\[
\tilde{Q}_{m_{1}}=:(\tilde{Q}(x), e^{i(2m_{1})\pi x})=\frac{c_{2m_{1}}}{i\pi(2m_{1})},\qquad 2m_{1}\neq0,
\]
\begin{equation}\label{d3}
\hat{Q}_{m_{1}}=:(\hat{Q}(x), e^{i(2m_{1}+1)\pi x})=\frac{c_{2m_{1}+1}}{i\pi(2m_{1}+1)},
\end{equation}
are denote the Fourier coefficients with respect to the systems  $\{e^{i(2m_{1})\pi x}:m_{1}\in \mathbb{Z}\}$ and $\{e^{i(2m_{1}+1)\pi x}:m_{1}\in \mathbb{Z}\}$ of the functions $\tilde{Q}(x)=\int_{0}^{x}\tilde{q}(t)\,dt$ and $\hat{Q}(x)=\int_{0}^{x}\hat{q}(t)\,dt$. Here, since $\tilde{Q}(1)=\int_{0}^{1}\frac{q(t)+q(1-t)}{2}\,dt$, we get $\tilde{Q}(1)=c_{0}=0$ if $q(1-t)=u$ use variable change. Similarly, $\hat{Q}(1)=0$. Now, using integration by parts, together by $\tilde{Q}(1)=0$ and $\hat{Q}(1)=0$, yields 
\[
\!\!\!\!\!\!\!\!\!\!\!\!\!\!\!\!\!\!\!\!\!\!\!\!\!\!\!\!\!\!\!\!\!\!\!\!\!\!\!b_{1}(\lambda_{m})=\frac{1}{2\pi^{2}m^{2}}\int_{0}^{1}(\tilde{q}^{2}(x)+\,(\tilde{Q}(x)-\tilde{Q}_{0})\,\tilde{q}^{\prime}(x))\,e^{i2m\pi x}\,dx
\]
\begin{equation}\label{d4}
+\frac{1}{2\pi^{2}m^{2}}\int_{0}^{1}(\hat{q}^{2}(x)+\,(\hat{Q}(x)-\hat{Q}_{0})\,\hat{q}^{\prime}(x))\,e^{i2m\pi x}\,dx+o(m^{-2}).
\end{equation}

Since $\tilde{q}(x)$ and $\hat{q}(x)$ are absolutely continuous a.e., $(\tilde{q}^{2}(x)+\,(\tilde{Q}(x)-\tilde{Q}_{0})\,\tilde{q}^{\prime}(x))\in L^{1}[0,1]$ and $(\hat{q}^{2}(x)+\,(\hat{Q}(x)-\hat{Q}_{0})\,\hat{q}^{\prime}(x))\in L^{1}[0,1]$. By using the Riemann-Lebesgue lemma, we obtain that 
\begin{equation}\label{d5}
b_{1}(\lambda_{m})=o(m^{-2}).
\end{equation} 
Now, by using (\ref{e9}) we prove that (see \cite{Kirac-inverse})
\begin{equation}\label{d7}
b_{2}(\lambda_{m})=o(m^{-2}).
\end{equation}
By considering that $\tilde{q}(x)$ and $\hat{q}(x)$ are absolutely continuous a.e., we obtain $c_{m_{1}}c_{m_{2}}c_{2m+m_{1}+m_{2}}=o(m^{-1})$ (see \cite{Shkalikov-Veliev}). This, together with (\ref{dif}) and (\ref{f7}), we write the sum on the right-hand-side of equation (\ref{e9}), as follows.
\[
\!\!\!\!\!\!\!\!\!\!\!\!\!\!\!\!\!\!\!\!\!\!\!\!\!\!\!\!\!\!\!\!\!\!\!\!\!\!\!\!\!\!\!\!\!\!\!\!b_{2}(\lambda_{m})=\sum_{m_{1},m_{2}=-\infty \atop m_{1},m_{1}+m_{2}\neq0,-2m}^{\infty}\frac{c_{m_{1}}c_{m_{2}}c_{2m+m_{1}+m_{2}}}{\prod_{t=1,2}[\lambda_{m}-(\pi(m+m_{t}))^{2}]}
\]
\[
=\frac{1}{\pi^{4}}\sum_{m_{1},m_{2}=-\infty \atop m_{1},m_{1}+m_{2}\neq0,-2m}^{\infty}\frac{c_{m_{1}}c_{m_{2}}c_{2m+m_{1}+m_{2}}}{-m_{1}(2m+m_{1})(-m_{1}-m_{2})(2m+m_{1}+m_{2})}+o(m^{-2})
\]
\[
|b_{2}(\lambda_{m})|=o(m^{-1})\sum_{m_{1},m_{2}}\frac{1}{|-m_{1}(2m+m_{1})(-m_{1}-m_{2})(2m+m_{1}+m_{2})|}
\]
\[
\!\!\!\!\!\!\!\!\!\!\!\!\!\!\!\!\!\!\!\!\!\!\!\!\!\!\!\!\!\!\!\!\!\!\!\!\!\!\!\!\!\!\!\!\!=o(m^{-1})O\left( \left( \frac{ln|m|}{m}\right) ^{2}\right)=o(m^{-2}).
\]
 Thus, the estimate of (\ref{d7}) is proved. 
\end{pf}

\begin{lem}\label{L3}
Let $q(x)\in W_{1}^{2}[0,1]$, $q^{(k)}(0)=q^{(k)}(1)$, $k=0,1$ and $c_{0}=0$. For all sufficiently large $m$, we obtain following estimates 
\begin{equation}\label{d19}
a_{1}(\lambda_{m})=\frac{1}{(2m)^{2}\pi^{2}}\int_{0}^{1}q^{2}(x)\,dx+o(m^{-2}),
\end{equation}
\begin{equation}\label{d8}
a_{2}(\lambda_{m})=o(m^{-2}).
\end{equation}
\end{lem}
\begin{pf}
In a similar way, arguing as in $b_{1}(\lambda_{m})$, let us prove that $a_{1}(\lambda_{m})$. For this, we write considering that (\ref{dif}) and (\ref{f7}) and the sum on the right-hand-side of equation (\ref{e6}) by using equalities (\ref{f1}) and (\ref{f2}), as following. 
\[
\!\!\!\!\!\!\!\!\!\!\!\!\!\!\!\!\!\!\!\!\!\!\!\!\!\!\!\!\!\!\!\!\!\!\!\!\!\!\!\!\!\!\!\!\!\!\!\!\!\!\!\!\!\!\!\!\!\!\!\!\!\!\!\!\!\!\!\!\!\!\!\!\!\!\!\!\!\!\!\!\!\!\!\!\!\!\! a_{1}(\lambda_{m})=\sum_{m_{1}=-\infty \atop m_{1}\neq0,-2m}^{\infty}\frac{c_{m_{1}}c_{m_{1}}}{[\lambda_{m}-(\pi(m+m_{1}))^{2}]}
\]
\[
\!\!\!\!\!\!\!\!\!\!\!\!\!\!\!\!\!\!\!\!\!\!\!\!\!\!\!\!\!\!\!\!\!\!\!\!\!\!\!\!\!\!\!\!\!\!\!\!\!\!\!\!\!\!\!\!\!\!\!\!\!\!\!\!\!\!=\frac{1}{\pi^{2}}\sum_{m_{1}=-\infty \atop m_{1}\neq0,-2m}^{\infty}\frac{c_{m_{1}}c_{m_{1}}}{-m_{1}(2m+m_{1})}+o(m^{-2})
\]
\[
\!\!\!\!\!\!\!\!\!\!\!\!\!\!\!\!\!\!\!\!\!\!\!\!\!\!\!\!\!\!\!\!\!\!\!\!\!\!\!\!\!\!\!\!\!\!\!\!\!\!\!\!\!\!\!\!\!\!\!\!\!\!\!\!\!\!\!\!\!\!\!\!\!\!\!\!\!\!\!\!\!
=\frac{1}{\pi^{2}}\sum_{2m_{1}\neq0,-2m}\frac{c_{2m_{1}}c_{2m_{1}}}{-2m_{1}(2m+2m_{1})}
\]
\begin{equation}\label{d9}
+\frac{1}{\pi^{2}}\sum_{2m_{1}+1\neq -2m}\frac{c_{2m_{1}+1}c_{2m_{1}+1}}{-(2m_{1}+1)(2m+2m_{1}+1)}+o(m^{-2}).
\end{equation}
It also follows from \cite[Lemma 2.2]{Kirac-inverse} (see \cite[Lemma 2.3]{Veliev}) that we get in our expression,
\[
\!\!\!\!\!\!\!\!\!\!\!\!\!\!\!\!\!\!\!\!\!\!\!\!\!\!\!\!\!\!\!\!\!\!\!\!\!\!\!\!\!\!\!\!\!\!\!\!\!\!\!\!\!\!\!\!\!\!\!\!\!\!\!\!\!\!\!\!\!\!\!\!\!\!\!\!\!\!\!\!\!\!\!\!\!\!\!\!\!\!\!\!\!\!\!\!\!\!\!\!\!\!\!\!\!\!\!\!\!\!\!\!\!\!\!\!\!\!\!\!\!\!\!\!\!\!\!\!\!\!\!\!\!\!\!\!\!\!\!\!\!\!\!\!\!\!\!\!\!\!\!\!\!\!\!\!\!\!\!\!\!\!\!\!\!\!\!\!\!\!\!\!\!\!\!\!\!\!\!\!a_{1}(\lambda_{m})
\]
\[
\!\!\!\!\!\!\!\!\!\!\!\!\!\!\!\!\!\!\!\!\!\!\!\!\!\!\!\!\!\!\!\!\!\!\!\!\!\!\!\!\!\!\!\!\!\!\!\!\!\!\!\!\!\!\!\!\!\!\!\!\!\!\!\!\!\!\!\!\!\!\!\!\!\!\!\!\!\!\!\!\!\!\!\!\!\!\!\!\!\!\!\!\!\!\!=\frac{2}{\pi^{2}}\sum_{2m_{1}\neq-2m, \atop 2m_{1}>0}\frac{c_{2m_{1}}c_{2m_{1}}}{(2m+2m_{1})(2m-2m_{1})}
\]
\[ 
\!\!\!\!\!\!\!\!\!\!\!\!\!\!\!\!\!\!\!\!\!\!\!\!\!\!\!\!\!\!\!\!\!\!\!\!\!\!+\frac{2}{\pi^{2}}\sum_{2m_{1}+1>0}\frac{c_{2m_{1}+1}c_{2m_{1}+1}}{(2m+2m_{1}+1)(2m-2m_{1}-1)}+o(m^{-2})
\]
\[
\!\!\!\!\!\!\!\!\!\!\!\!\!\!\!\!\!\!\!\!\!\!\!\!\!\!\!\!\!\!\!\!\!\!\!\!\!\!\!\!\!\!\!\!\!\!\!\!\!\!\!\!\!\!\!\!\!\!\!\!\!\!\!\!\!\!\!\!\!\!\!\!\!\!\!\!\!\!\!\!\!\!\!\!=-\int_{0}^{1}(\tilde{G}^{+}(x,m)-\tilde{G}_{0}^{+}(m))^{2}\,e^{i(-4m)\pi x}\,dx
\]
\[
 \!\!\!\!\!\!\!\!\!\!\!\!\!\!\!\!\!\!\!\!\!\!\!\!\!\!-\int_{0}^{1}(\hat{G}^{+}(x,m)-\hat{G}_{0}^{+}(m)+o(m^{-2}))^{2}\,e^{i(-4m)\pi x}\,dx+o(m^{-2})
\]
\[
\!\!\!\!\!\!\!\!\!\!\!\!\!\!\!\!\!\!\!\!\!\!\!\!\!\!\!\!\!\!\!\!\!\!\!\!\!\!\!\!\!\!\!\!\!\!\!\!\!\!\!\!\!\!\!\!\!\!\!\!\!\!\!\!\!\!\!\!\!\!\!\!\!\!\!\!\!\!\!\!\!=-\int_{0}^{1}(\tilde{G}^{+}(x,m)-\tilde{G}_{0}^{+}(m))^{2}\,e^{i(-4m)\pi x}\,dx
\]
\[
\!\!\!\!\!\!\!\!\!\!\!\!\!\!\!\!\!\!\!\!\!\!\!\!\!\!\!\!\!\!\!\!\!\!\!\!\!\!\!\!\!\!\!\!\!\!\!\!\!\!\!-\int_{0}^{1}(\hat{G}^{+}(x,m)-\hat{G}_{0}^{+}(m))^{2}\,e^{i(-4m)\pi x}\,dx+o(m^{-2})
\]
\[
=\frac{1}{i\pi(-4m)}\int_{0}^{1}2(\tilde{G}^{+}(x,m)-\tilde{G}_{0}^{+}(m))\,(\tilde{q}(x)\,e^{-i(-2m)\pi x}-c_{(-2m)})\,e^{i(-4m)\pi x}\,dx
\]
\begin{equation}\label{d10}
+\frac{1}{i\pi(-4m)}\int_{0}^{1}2(\hat{G}^{+}(x,m)-\hat{G}_{0}^{+}(m))\,\hat{q}(x)\,e^{-i(-2m)\pi x}\,e^{i(-4m)\pi x}\,dx +o(m^{-2})
\end{equation}
where
\[
\tilde{G}^{\pm}_{m_{1}}(m)=:(\tilde{G}^{\pm}(x,m), e^{i(2m_{1})\pi x})=\frac{c_{2m_{1}\pm(-2m)}}{i\pi(2m_{1})},
\]
\begin{equation}\label{d11}
\hat{G}^{\pm}_{m_{1}}(m)=:(\hat{G}^{\pm}(x,m), e^{i(2m_{1}+1)\pi x})=\frac{c_{2m_{1}+1\pm(-2m)}}{i\pi(2m_{1}+1)}+\frac{2}{(2m_{1}+1)^{2} \pi ^{2}}\int_{0}^{1}\hat{q}(t)e^{i(2m) \pi t}\,dt
\end{equation}
for $2m_{1}\neq0$ are denote the Fourier coefficients with respect to system $\{e^{i(2m_{1})\pi x}:m_{1}\in \mathbb{Z}\}$ and  $\{e^{i(2m_{1}+1)\pi x}:m_{1}\in \mathbb{Z}\}$ of the functions 
\[
\tilde{G}^\pm(x,m)=\int_{0}^{x}\tilde{q}(t)\,e^{\mp i(-2m)\pi t}\,dt-c_{\pm(-2m)}x , 
\]
\begin{equation}\label{d12}
\hat{G}^\pm(x,m)=\int_{0}^{x}\hat{q}(t)\,e^{\mp i(-2m)\pi t}\,dt-x\int_{0}^{1} \hat{q}(t)\,e^{\mp i(-2m)\pi t}\,dt
\end{equation}
and
\[
\tilde{G}^{\pm}(x,m)-\tilde{G}_{0}^{\pm}(m)=\sum_{2m_{1}\neq-2m}\frac{c_{2m_{1}}}{i\pi(2m_{1}\mp(-2m))}e^{i(2m_{1}\mp(-2m))\pi x} ,
\]
\begin{equation}\label{d25}
\hat{G}^{\pm}(x,m)-\hat{G}_{0}^{\pm}(m)=\sum \frac{c_{2m_{1}+1}}{i\pi(2m_{1}+1 \mp(-2m))}e^{i(2m_{1}+1 \mp(-2m))\pi x}+o(m^{-2}).
\end{equation}

Here, in the second expression of (\ref{d25}), the second term on the right-hand side of the equality is uniform in $x$. Considering the \cite[Lemma 1]{Kirac-Titchmarsh} and (\ref{d12}) we have following the estimates 
\[
\tilde{G}^{\pm}(x,m)-\tilde{G}_{0}^{\pm}(m)=\tilde{G}^{\pm}(x,m)-\int_{0}^{1}\tilde{G}^{\pm}(x,m)\,dx=o(1),
\]
\begin{equation}\label{d13}
\hat{G}^{\pm}(x,m)-\hat{G}_{0}^{\pm}(m)=\hat{G}^{\pm}(x,m)-\int_{0}^{1}\hat{G}^{\pm}(x,m)\,dx=o(1),\quad \text{as} \qquad m\rightarrow\infty
\end{equation}
uniformly in $x$.

Taking into account the equalities (see (\ref{d12}))
\begin{equation}\label{d14}
\tilde{G}^{\pm}(1,m)=\tilde{G}^{\pm}(0,m)=0,\qquad \hat{G}^{\pm}(1,m)=\hat{G}^{\pm}(0,m)=0,
\end{equation}
and since $\tilde{q}(x)$ and $\hat{q}(x)$ are absolutely continuous a.e., using the integration by parts obtain for the right-hand side of (\ref{d10}), the value 
\[
\!\!\!\!\!\!\!\!\!\!\!\!\!\!\!\!\!\!\!\!\!\!\!\!\!\!\!\!\!\!\!\!\!\!\!\!\!\!\!\!\!\!\!\!\!\!\!\!\!\!\!\!\!\!\!\!\!\!\!\!\!\!\!\!\!\!\!\!\!\!\!\!\!\!\!\!\!\!\!\!\!
\!\!\!\!\!\!\!\!\!\!\!\!\!\!\!\!\!\!\!\!\!\!\!\!\!\!\!\!\!\!\!\!\!\!\!\!\!\!\!\!\!\!\!\!\!\!\!\!\!\!\!\!\!\!
\!\!\!\!\!\!\!\!\!\!\!\!\!\!\!\!\!\!\!\!\!\!\!\!\!\!\!\!\!\!\!\!\!\!\!\!\!\!a_{1}(\lambda_{m})
\]
\[
=\frac{1}{(2m)^{2}\pi^{2}}\left[ \int_{0}^{1}\tilde{q}^{2}+\int_{0}^{1}(\tilde{G}^{+}(x,m)-\tilde{G}_{0}^{+}(m))\,\tilde{q}^{\prime}(x)\,e^{i(-2m)\pi x}\,dx\right]-\frac{3|c_{(-2m)}|^{2}}{2\pi^{2}(2m)^{2}}
\]
\[
+\frac{1}{(2m)^{2}\pi^{2}}\left[ \int_{0}^{1}\hat{q}^{2}+
\int_{0}^{1}(\hat{G}^{+}(x,m)-\hat{G}_{0}^{+}(m))\,\hat{q}^{\prime}(x)\,e^{i(-2m)\pi x}\,dx\right]+o(m^{-2})
\]
for sufficiently large $m$. Thus, by using the Riemann - Lebesgue lemma, together with $(\tilde{G}^{+}(x,m)-\tilde{G}_{0}^{+}(m))\,\tilde{q}^{\prime}(x)\in L^{1}[0,1]$ and $(\hat{G}^{+}(x,m)-\hat{G}_{0}^{+}(m))\,\hat{q}^{\prime}(x)\in L^{1}[0,1]$, we get (\ref{d19}).

Now, let's prove that $a_{2}(\lambda_{m})=o(m^{-2})$. Similarly, together with (\ref{e7}), (\ref{dif}) and (\ref{f7}), we get 
\[
a_{2}(\lambda_{m})=\sum_{ m_{1},m_{2}=-\infty  \atop m_{1},m_{1}+m_{2}\neq0,-2m}^{\infty}\frac{c_{m_{1}}c_{m_{2}}c_{m_{1}+m_{2}}}{\prod_{t=1,2}[\lambda_{m}-(\pi(m+m_{t}))^{2}]}
\]
\begin{equation}\label{d15}
=\sum_{m_{1},m_{2}}\frac{\pi^{-4}c_{m_{1}}c_{m_{2}}c_{m_{1}+m_{2}}}{-m_{1}(2m+m_{1})(-m_{1}-m_{2})(2m+m_{1}+m_{2})}+o(m^{-2}).
\end{equation}
Arguing as in \cite[Lemma 4]{Kirac-Titchmarsh} and \cite{Kirac-inverse}, using the summation varyant $m_{2}$ to impression the previous $m_{1}+m_{2}$ in (\ref{d15}), we write (\ref{d15}) in the formula
\[
 a_{2}(\lambda_{m})=\frac{1}{\pi^{4}}\sum_{m_{1},m_{2}}\frac{c_{m_{1}}c_{m_{2}-m_{1}}c_{m_{2}}}{-m_{1}(2m+m_{1})(-m_{2})(2m+m_{2})}+o(m^{-2}).
 \]
Here, the forbidden indices in the sums take the form of $m_{1},m_{2}\neq0,-2m$. By the equality 
 \[
 \frac{1}{-k(2m+k)}=\frac{1}{2m}\left( \frac{1}{2m+k}-\frac{1}{k}\right) 
 \]
 we have 
\begin{equation}\label{d16}
 a_{2}(\lambda_{m})=\frac{1}{\pi^{4}(2m)^{2}}\sum_{j=1}^{4}S_{j},
\end{equation}
where 
\[
S_{1}=\sum_{m_{1},m_{2}}\frac{c_{m_{1}}c_{m_{2}-m_{1}}c_{m_{2}}}{m_{1}m_{2}},\qquad S_{2}=\sum_{m_{1},m_{2}}\frac{c_{m_{1}}c_{m_{2}-m_{1}}c_{m_{2}}}{-m_{2}(2m+m_{1})}
\]
\[
S_{3}=\sum_{m_{1},m_{2}}\frac{c_{m_{1}}c_{m_{2}-m_{1}}c_{m_{2}}}{-m_{1}(2m+m_{2})},\qquad S_{4}=\sum_{m_{1},m_{2}}\frac{c_{m_{1}}c_{m_{2}-m_{1}}c_{m_{2}}}{(2m+m_{1})(2m+m_{2})}.
\]
Now, using the first equality in (\ref{d3}), integration by parts and the assumption $c_{0}=0$ which means $\tilde{Q}(1)=0$, we get for even that
\[
\tilde{S}_{1}=\pi^{2}\int_{0}^{1}(\tilde{Q}(x)-\tilde{Q}_{0})^{2}\,\tilde{q}(x)\,dx=0.
\]
Similarly, taking into account (\ref{d3}) and (\ref{d11})-(\ref{d14}), we obtain for even by the Riemann-Lebesgue lemma, the following relations
\[
\tilde{S}_{2}=-\pi^{2}\int_{0}^{1}(\tilde{Q}(x)-\tilde{Q}_{0})\,(\tilde{G}^{+}(x,m)-\tilde{G}_{0}^{+}(m))\,\tilde{q}(x)\,e^{i\pi(-2m)x}\,dx=o(1),
\]
\[
\tilde{S}_{3}=-\pi^{2}\int_{0}^{1}(\tilde{Q}(x)-\tilde{Q}_{0})\,(\tilde{G}^{-}(x,m)-\tilde{G}_{0}^{-}(m))\,\tilde{q}(x)\,e^{i\pi(2m)x}\,dx=o(1)
\]
and with the first equality in (\ref{d13}), 
\[
\tilde{S}_{4}=\pi^{2}\int_{0}^{1}(\tilde{G}^{-}(x,m)-\tilde{G}_{0}^{-}(m))\,(\tilde{G}^{+}(x,m)-\tilde{G}_{0}^{+}(m))\,\tilde{q}(x)\,dx=o(1).
\]

In the same way, using the second equality in (\ref{d3}), integration by parts and the assumption $c_{0}=0$ which means $\hat{Q}(1)=0$, we obtain for odd that
\[
\hat{S}_{1}= \pi^{2} \int_{0}^{1}(\hat{Q}(x)-\hat{Q}_{0})^{2}\hat{q}(x)\,dx=0.
\]

Similarly, considering (\ref{d3}) and (\ref{d11})-(\ref{d14}), we get for odd by the Riemann-Lebesgue lemma, the following relations
\[
\hat{S}_{2}=-\pi^{2}\int_{0}^{1}(\hat{Q}(x)-\hat{Q}_{0})\,(\hat{G}^{+}(x,m)-\hat{G}_{0}^{+}(m))\,\hat{q}(x)\,e^{i\pi(-2m)x}\,dx=o(m^{-2}),
\]
\[
\hat{S}_{3}=-\pi^{2}\int_{0}^{1}(\hat{Q}(x)-\hat{Q}_{0})\,(\hat{G}^{-}(x,m)-\hat{G}_{0}^{-}(m))\,\hat{q}(x)\,e^{i\pi(2m)x}\,dx=o(m^{-2})
\]
and with the second equality in (\ref{d13}), 
\[
\hat{S}_{4}=\pi^{2}\int_{0}^{1}(\hat{G}^{-}(x,m)-\hat{G}_{0}^{-}(m))\,(\hat{G}^{+}(x,m)-\hat{G}_{0}^{+}(m))\,\hat{q}(x)\,dx=o(m^{-2}).
\] 
Thus, (\ref{d16}) proved.
\end{pf}

\renewenvironment{pf}[1]{{\bfseries Proof of Theorem \ref{maintheorem}.#1}}{}
\begin{pf}
First, considering (\ref{mc0}) we find $c_{0}=0$ since \[
\lambda_{m}=(m \pi)^{2}+c_{0}+o(1).
\]
By substituting the estimates obtained in $c_{0}=0$,  $c_{2m}=o(m^{-2})$, Lemma \ref{L2} and Lemma \ref{L3} in Lemma \ref{L1}, we obtain the form
\[
\lambda_{m}=(m \pi)^{2}+\frac{1}{(2m)^{2}\pi^{2}}\int_{0}^{1}q^{2}(x) \,dx+o(m^{-2}).
\]
From the hypothesis, 
\[
\int_{0}^{1}q^{2}(x)\,dx=0,
\]
so  $q=0$ a.e. 
\end{pf}

\section*{Acknowledgments}
The author is grateful to the anonymous reviewers for their helpful
comments and suggestions.

\bibliographystyle{elsarticle-num}
\bibliography{article}


\end{document}